\documentclass[runningheads,a4paper]{llncs}

\usepackage{amssymb,amsfonts}
\usepackage{graphicx}

\usepackage{url}
\newcommand{\keywords}[1]{\par\addvspace\baselineskip
\noindent\keywordname\enspace\ignorespaces#1}

\begin{document}

\mainmatter  % start of an individual contribution

\title{Sub-Riemannian Fast Marching in SE(2)}

\titlerunning{Sub-Riemannian Fast Marching in SE(2)}

\author{G. Sanguinetti$^1$\thanks{The research leading to the results of this article has received funding from the European Research Council under the (FP7/2007-€"2014-)ERC grant ag.~no.~335555 and from (FP7-PEOPLE-2013-ITN)EU Marie-Curie ag.~no.~607643.}, E. Bekkers$^2$, R. Duits$^{1,2}$, M.H.J. Janssen$^1$, A. Mashtakov$^2$, J.M. Mirebeau$^3$}

%\author{Anonymous}

%
\authorrunning{Sanguinetti et al.}

\institute{ Eindhoven University of Technology, The Netherlands, \\
$^1$Dept. of Mathematics and Computer Science, $^2$Dept. of Biomedical Engineering. \\ 
 $^3$University Paris-Dauphine, Laboratory Ceremade, CNRS
\email{\{G.R.Sanguinetti,E.J.Bekkers,R.Duits,M.H.J.Janssen,A.Mashtakov\}@tue.nl}\\
 \email{mirebeau@ceremade.dauphine.fr}}

%\institute{ Institute, \email{email}}

\maketitle
\vspace{-1.5ex}
\begin{abstract}
We propose a Fast Marching based implementation for computing sub-Riemanninan (SR) geodesics in the roto-translation group SE(2), with a metric depending on a cost induced by the image data. The key ingredient is a Riemannian approximation of the SR-metric. Then, a state of the art Fast Marching solver that is able to deal with extreme anisotropies is used to compute a SR-distance map as the solution of a corresponding eikonal equation. Subsequent backtracking on the distance map gives the geodesics. To validate the method, we consider the uniform cost case in which exact formulas for SR-geodesics are known and we show remarkable accuracy of the numerically computed SR-spheres. We also show a dramatic decrease in computational  time with respect to a previous PDE-based iterative approach. Regarding image analysis applications, we show the potential of considering these data adaptive geodesics for a fully automated retinal vessel tree segmentation.
\end{abstract}
\section{Introduction}
In this article we study a curve optimization problem in the space of coupled positions and orientations $\mathbb{R}^2\!\times\! S^1$, which we identify with  roto-translation group $SE(2)$. % in the roto-translation group endowed with sub-Riemannian metric which depends on an external given cost. 
We aim to compute the shortest curve $\gamma(t)=(x(t),y(t),\theta(t))\in SE(2)$  that connects 2 points $\gamma(0)=(x_0,y_0,\theta_0)$ and $\gamma(L)=(x_1,y_1,\theta_1)$  with the restriction that the curve is "lifted" from a planar curve in the sense that the third variable $\theta$ is given by $\theta(t)=\arg(\dot{x}(t)+i\, \dot{y}(t))$, see Fig. 1. This restriction imposes a so-called sub-Riemannian (through the text we denote sub-Riemannian as SR) metric that constrains the tangent vectors to the curve  to be contained in a subspace of the tangent space at every point. This subspace is a plane, which differs from point to point, and is the set of all possible tangents to curves in $SE(2)$ that are lifted from smooth planar curves. A SR-metric is a degenerated Riemannian metric in which one direction, the one perpendicular to the path in this case, is prohibited (i.e. it has infinite cost). On top of the SR-metric, we also consider a smooth external cost, which weights the metric tensor in every location and allows for data-adaptivity. 
\begin{figure}
\centering
\includegraphics[width=.34\textwidth]{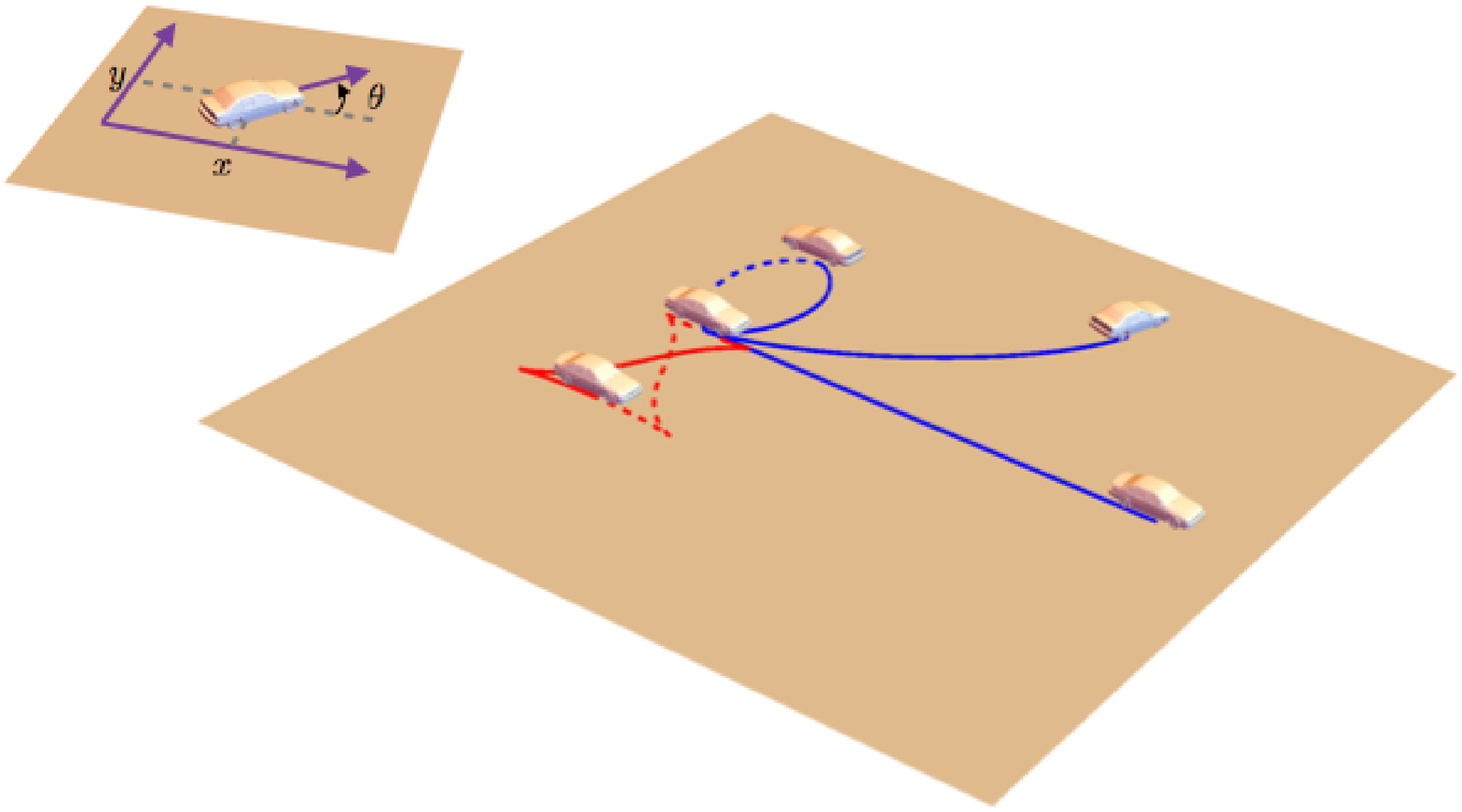}
\includegraphics[width=.3\textwidth]{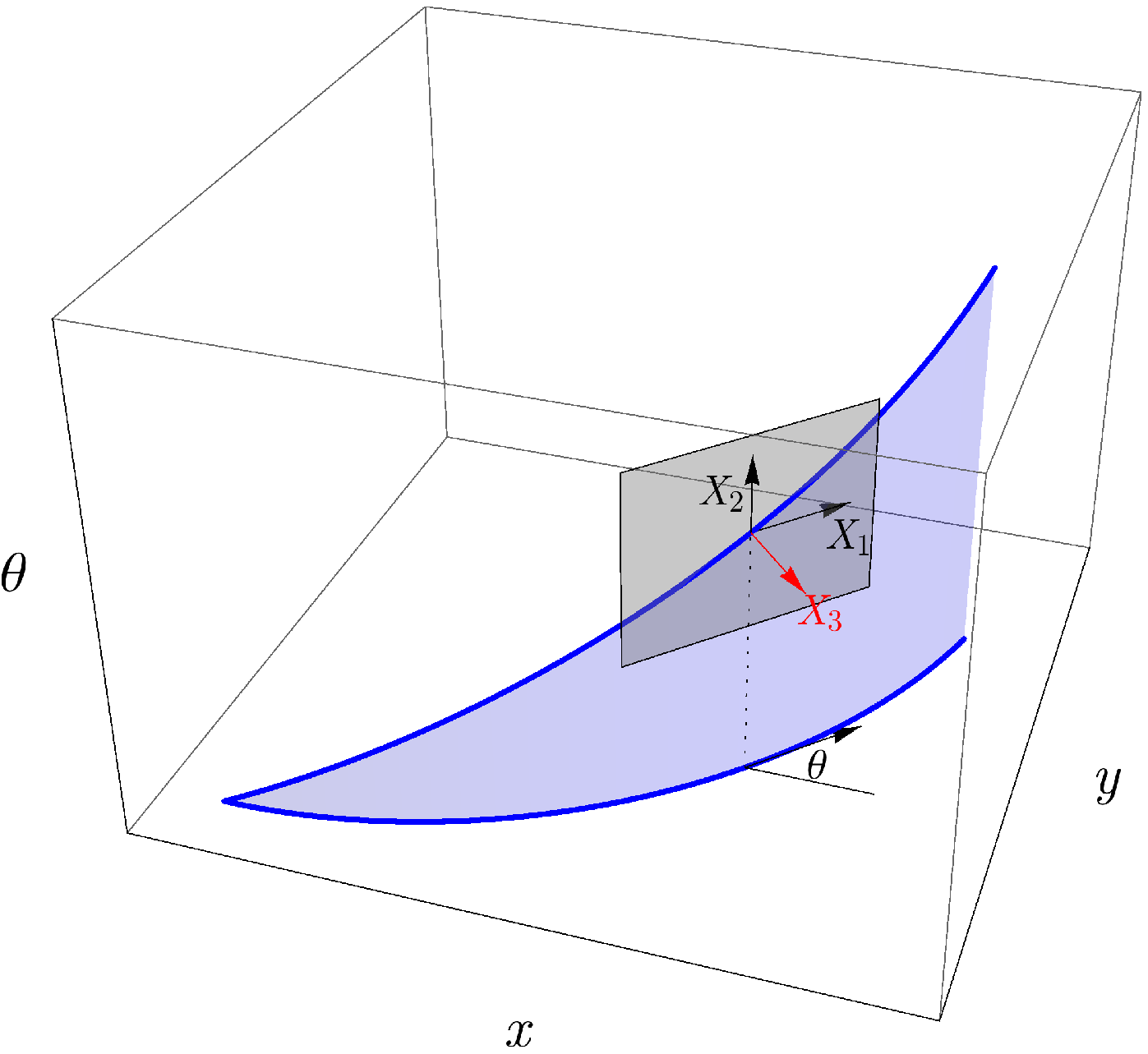}
\includegraphics[width=.34\textwidth]{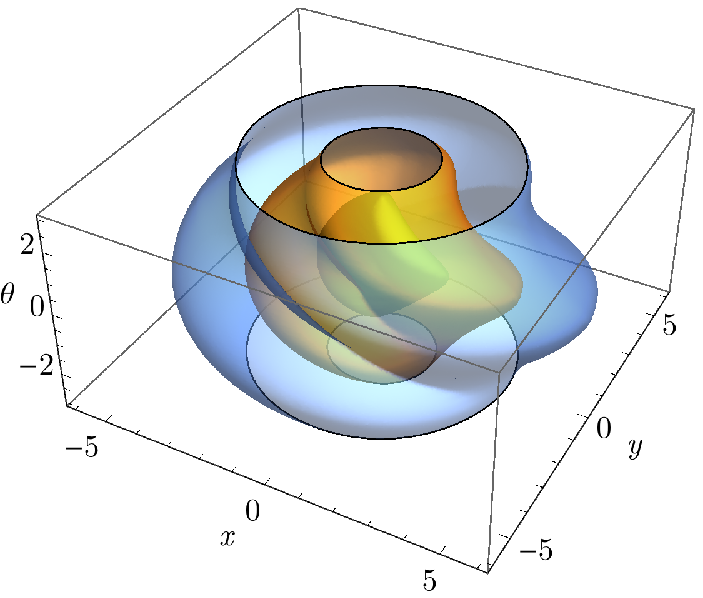}
\caption{\textbf{Left:} the sub-Riemannian problem in $SE(2)$ can be identified with that of a car with two controls (giving gas and steering the wheel). \textbf{Center:} the paths are ''lifted'' into curves in $SE(2)=\mathbb{R}^2\!\times\! S^1$ with tangent vectors constrained to the plane spanned by the vector fields $X_1$ and $X_2$ (eq.~(\ref{leftinvariant})) associated with the controls. \textbf{Right:} the SR-spheres (for $t=2,4$ and $6$) obtained via the FM-LBR method.}
\label{fig:Fig1}
\end{figure}

Essentially, the SR-problem in $SE(2)$ is that of a car that can go forward, backward and rotate (a so-called Reeds-Shepp car% \cite{reeds}
) so the possible states of the car form a 3D manifold given by the position $(x,y)$ and the orientation $\theta$ of the car. Then, admissible trajectories of the car are parametrized by only 2 control variables associated to the car moving along a straight line (giving gas) and to a change of direction (turning the steering wheel). The fact that the car cannot step aside infinitesimally imposes the SR-geometry. %This can be generalized by introducing a cost function on the state space. 
Finally, the curve optimization problem is to find among all possible trajectory between two given states, the one with  minimal SR-length.

In image analysis the SR-geodesics in $SE(2)$ have been proposed in \cite{cittisarti} as candidates for completion of occluded contours. Here, the geometrical structure is used as a model for the functional architecture of the primary visual cortex, extending the work in \cite{petitot} to the $SE(2)$ group.  This model has proven to be valuable in numerous applications \cite{benyosef,cittisarti,duits,franken}, and it becomes powerful when combined with the orientation score theory \cite{franken} that allows for an invertible stable transformation between 2D images and functions on the SE(2) group.
The main advantage of considering this space of positions and orientations is that intersecting curves are automatically disentangled, and therefore the processing in the extended domain naturally deals with complex structures such as crossing.

Sub-Riemannian geodesics in the uniform cost case (the same cost for all the $SE(2)$ elements) were studied by several authors (e.g. \cite{duits,yuri}). %Usually, the boundary value problem is solved via shooting methods. 
Recently, a wavefront based method for computing SR-geodesic that also deals with the non-uniform cost case has been proposed in \cite{bekkers}. This method is an extension to the SR-case of a classical PDE-approach for computing cost-adaptive geodesics used in computer vision, where the metric tensor is induced by the image itself. The main idea is to consider propagation of equidistant surfaces described by level sets of the viscosity solution of an eikonal equation, while subsequent backtracking gives the geodesics. 
%In \cite{bekkers} both the SR-eikonal equation and the backtracking are derived via Pontryagin's Maximum Principle. Also in \cite{bekkers}, the authors rely on a computationally expensive iterative approach for computing the solution to the eikonal equation based on left-invariant finite-difference discretization of the PDE combined with a suitable upwind scheme.
In order to solve the eikonal equation, the authors rely on a computationally expensive iterative approach based on a left-invariant finite-difference discretization of the PDE combined with a suitable upwind scheme. 

Here, again in the spirit of computer vision methods, we aim to compute the SR-distance map using a Fast Marching method \cite{sethian}. This technique, closely related to Dijkstra's method for computing the shortest paths on networks, allows for a significant speed up in the computation of the eikonal equation's viscosity solution. The main difficulty in our case is that classical solvers are unable to deal with the extreme (degenerated) anisotropy of the SR-metric. Recently, a modification of the Fast Marching method using lattice basis reduction (FM-LBR) that solves this problem was introduced in \cite{mirebeau} (code available at \url{https://github.com/Mirebeau/ITKFM}).  Then, the purpose of this paper is to show how the SR-curve optimization problem can be numerically solved  using the FM-LBR method. The key aspects to consider are a Riemannian relaxation of the SR-problem and expressing the resulting metric tensor as a matrix-induced Riemannian metric in a fixed Cartesian frame.  We develop these ideas in the Theory section. Then, two experiments are presented. The first one considers the uniform cost case ($\mathcal{C}\!=1\!$) and shows that the FM-LBR based method presented here outperforms the iterative implementation in \cite{bekkers} in terms of CPU time, while keeping a similar accuracy. The advantages of considering data-adaptive SR-geodesics for extracting blood vessels in retinal images are illustrated in the second experiment. 
%
%... TEXT NEEDED HERE
%This article is organized as follows, Section \ref{theory} describe the precise problem formulation, the sub-Riemannian eikonal equation and how ... TEXT NEEDED HERE
%We aim to compute the curve with minimal length that connect 2 point  with the restriction that the curve is the "lifting" of a planar curve. With lifting 
%A typical example of such a curve is the path described by a car 
\vspace{-1.5ex}
\section{Theory} \label{theory}\vspace{-1.5ex}
\textbf{Problem formulation.} Let $g=(x,y,\theta)$ be an element of $SE(2)=\mathbb{R}^2\!\rtimes \!S^1$. A natural moving frame of reference in $SE(2)$ is described by the left-invariant vector fields  $\{X_{1},X_{2},X_{3}\}$ spanning the tangent space at each element $g$:
\begin{equation} \label{leftinvariant}
\begin{array}{lll}
 X_{1}= \cos \theta \partial_{x} +\sin \theta
 \partial_{y},& 
 X_{2}= \partial_{\theta},& 
 X_{3}= -\sin \theta \partial_{x}+\cos \theta
 \partial_{y}.
\end{array}
\end{equation}
The tangents $\dot \gamma (t)$ along curves $\gamma(t)=(x(t),y(t),\theta(t))\in SE(2)$ can be written as $\dot \gamma (t)=\sum_{i=1}^3 u^i(t)\left.X_i\right|_{\gamma(t)}$. Only the curves with $u^3=0$ are liftings of planar curves (see fig. \ref{fig:Fig1}). Then, the tangents to curves that are liftings of planar curves are expressed as $\dot \gamma(t)=u^1(t)\left.X_1\right|_{\gamma(t)}+u^2(t)\left.X_2\right|_{\gamma(t)}$ and they span the so-called distribution $\Delta=\mathrm{span}\left\{X_1,X_2\right\}$. Now, the triplet $\left(SE(2),\Delta,G_{0}\right)$ defines a sub-Riemannian manifold with inner product $G_{0}$ given by:
\begin{equation} \label{eq:SRmetric}
G_{0}(\dot{\gamma},\dot{\gamma}) = \mathcal{C}(\gamma)^2\left(\beta^2 |\dot{x}\cos\theta  \!+\!\dot{y}\sin\theta |^2 + |\dot{\theta}|^2\right).
\end{equation}
In view of the car example (see Fig. \ref{fig:Fig1}), the parameter $\beta>0$ balances between the cost of giving gas (the $X_1$ direction) and the turning of the wheel (the $X_2$ direction). A smooth function $\mathcal{C}:SE(2)\to [\delta,1], \, \delta>0$ is the given external cost. In order to keep the notation sober during this paper, we do not indicate the dependence of $G_0$ on the cost $\mathcal{C}$ nor that it also depends on the curve $\gamma$. Note that the subindex $0$ in the metric tensor recalls the SR-structure imposed by allowing  null displacement in the $X_3$ direction (i.e. infinite cost for the car to move aside). This choice of notation will become clear later the text. Now, optimal paths $\gamma$ of the car in our extended position orientation space are solutions of the following problem:
\begin{equation}\label{eq:min}
W(g)=d_{0}
(g,e)=\min\left\{\left.\int\limits_0^T \| \dot\gamma(t) \|_{0}\,\mathrm{ dt} \,\right |\, \dot\gamma\in\Delta, \gamma(0)=e,\gamma(T)=g,T\geq0 \right\}
\end{equation}
where $e=(0,0,0)$ is the origin, where
%\begin{equation}
$\left\|\dot\gamma(t) \right\|_{0}=\sqrt{G_{0}(\dot\gamma(t),\dot\gamma(t))}$
%\end{equation}
 is the SR-norm and $d_{0}$ is the SR-distance on $SE(2)$.\\

\textbf{Riemannian approximation.} It is possible to obtain a Riemannian approximation of the SR-problem by expanding the metric tensor in eq. (\ref{eq:SRmetric}) to:
\begin{equation} \label{eq:metric}
G_{\epsilon}(\dot{\gamma},\dot{\gamma}) = G_{0}(\dot{\gamma},\dot{\gamma}) + \epsilon^{-2}\mathcal{C}(\gamma)^2\beta^2 |\dot{x}\sin\theta  \!-\!\dot{y}\cos\theta |^2,
\end{equation}
where $\epsilon$ determines the amount of anisotropy between $X_3$ and $\Delta$. This definition bridges the SR-case, obtained at the limit $\epsilon\downarrow 0$, with the full Riemannian metric tensor when $\epsilon=1$ (isotropic in the spatial directions $X_1$ and $X_3$). Actually, it is easy to verify that if $\mathcal{C}=1$ and $\beta=\epsilon=1$, then
$G_{1}(\dot{\gamma},\dot{\gamma}) =  |\dot{x}|^2 +   |\dot{y}|^2 + |\dot{\theta}|^2.$
Also, by replacing $G_0$ with $G_\epsilon$ in eq.~(\ref{eq:min}) we can construct a $SE(2)$ Riemannian norm $\|\cdot\|_\epsilon$ and a $SE(2)$ Riemannian distance $d_\epsilon$ satisfying $\lim\limits_{\epsilon \downarrow 0}\|\cdot\|_\epsilon = \|\cdot\|_0$ and $\lim\limits_{\epsilon \downarrow 0}d_\epsilon = d_0$. 
\begin{figure}
\centering
\includegraphics[width=.25\textwidth]{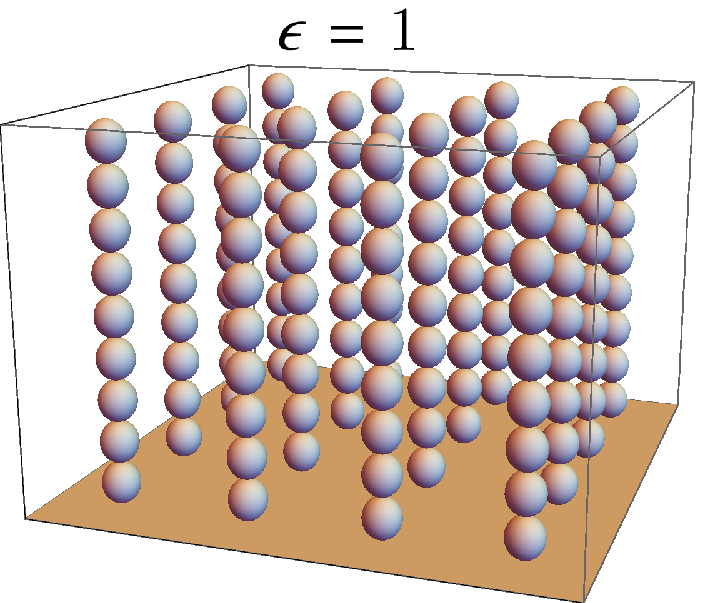}
\includegraphics[width=.25\textwidth]{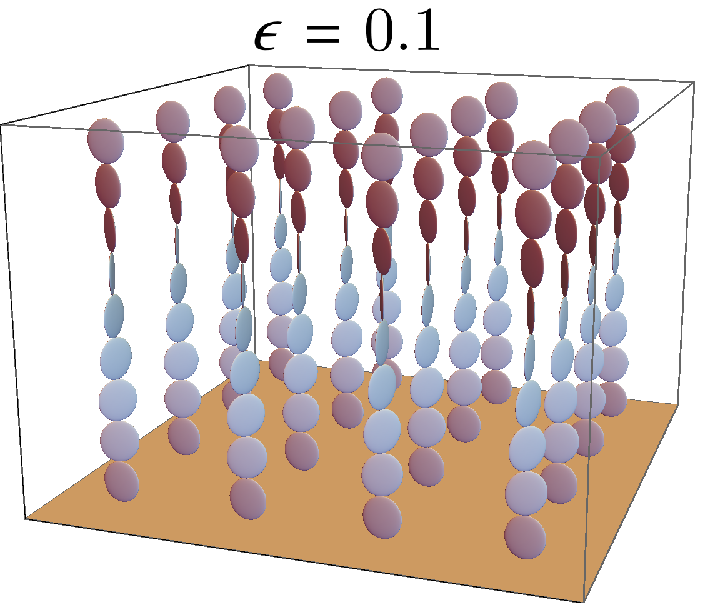}
\includegraphics[width=.25\textwidth]{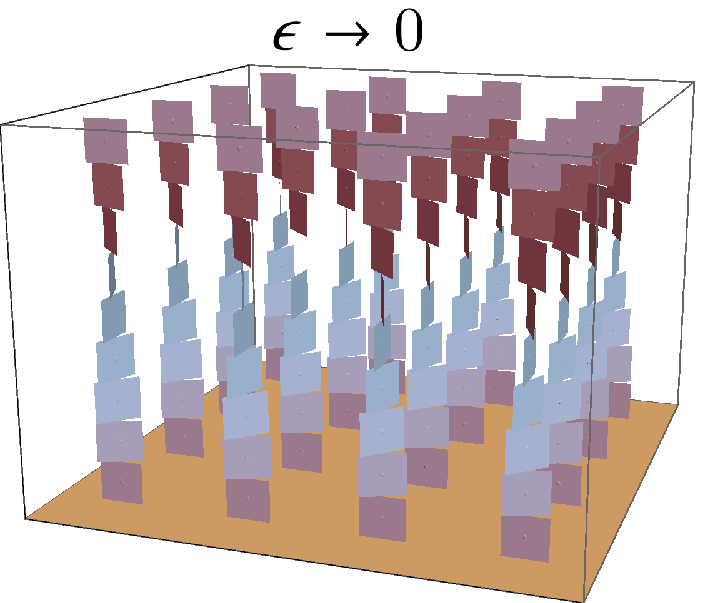}
\caption{ Each ellipsoid represents the Tissot's indicatrix of the metric $G_\epsilon$ at different elements $g\in SE(2)$ (for the case $\mathcal{C}(g)\!=\!1$ and $\beta\!=\!1$). The parameter $\epsilon$ in eq.~(\ref{eq:metric}) bridges the Riemannian case with the SR-one.  When $\epsilon=1$ each direction has the same cost. At the limit $\epsilon \downarrow 0$, the direction $X_3$ has infinite cost and the distribution $\Delta$ appears.}
\label{fig:metricapprox}
\end{figure}
%
%
%
%In this Riemannian case, one can find a matrix representation $\mathcal{M}_\epsilon^\mathcal{C}$ of the metric $G_\epsilon^\mathcal{C}$ given by:
%\begin{equation}\mathcal{M}_\epsilon^\mathcal{C}=\frac{1}{\mathcal{C}^2}
%\left(\begin{array}{ccc}
%\beta^2(\cos^2\theta + \epsilon^{-2}\sin^2\theta)& \beta^2(1- \epsilon^{-2})\cos^2\theta \sin\theta& 0\\
%\beta^2(1- \epsilon^{-2})\cos^2\theta \sin\theta & \beta^2(\epsilon^{-2} \cos^2\theta + \sin^2\theta) & 0\\
%0 & 0 & 1
%\end{array}\right).
%\end{equation}
% The set of tangents to these curves form the so-called horizontal distribution $\Delta=\textrm{span}\{X_{1},X_{2}\}$.  
%The inner product $G_0^\mathcal{C}: SE(2) \times \Delta \times \Delta \to \mathbb{R}$ is given by
%\begin{equation} \label{metric}
%G_0^\mathcal{C}(\dot{\gamma}(t),\dot{\gamma}(t)) = \beta^2 |\dot{x}(t)\cos\theta(t)  \!+\!\dot{y}(t)\sin\theta(t) |^2 + |\dot{\theta}(t)|^2,
%\end{equation}
%where $\beta>0$ constant set the balance between the spatial and angular components, and the smooth function $\mathcal{C}:SE(2)\to [\delta,1], \, \delta>0$ is the given external cost. Then, the triplet .
%In view of the car example (see Fig. \ref{fig:Fig1}), the parameter $\beta$ balances between the cost of giving gas (the $X_1$ direction) and the turning of the wheel (the $X_2$ direction). Now the optimal path of the car in our extended position orientation space is the minimizer of the of the following problem:

\textbf{Solution via the eikonal equation.} Now we can present the eikonal system that solves the problem~(\ref{eq:min}) by computing the distance map $W(g)$ as proved in \cite{bekkers}. 
Following \cite{cittisarti}, let us introduce some differential operators that will simplify the notation of the remaining equations.
Let $U$ be a smooth function $U:SE(2)\to \mathbb{R}$. The Riemannian %(SR at the limit $\epsilon\to0$) 
gradient $\nabla^\epsilon$, computed as the inverse of the metric tensor $G_\epsilon$ acting on the derivative, is given by:
\begin{equation}
\nabla^\epsilon U := G_\epsilon^{-1}\mathrm{d}U= \mathcal{C}^{-2}\beta^{-2}(X_1 U)X_1 +\mathcal{C}^{-2}(X_2 U)X_2+\epsilon^2\mathcal{C}^{-2}\beta^{-2}(X_3 U)X_3.
\end{equation}
Then, the norm of the gradient $\nabla^\epsilon$ is given by:
\begin{equation}
\left\|\nabla^\epsilon U\right\|_{\epsilon} = \sqrt{\mathcal{C}^{-2}\beta^{-2}|X_1 U|^2 +\mathcal{C}^{-2}|X_2 U|^2+\epsilon^2\mathcal{C}^{-2}\beta^{-2}|X_3U|^2}.
\end{equation}
Thus, the eikonal system that characterizes the propagation of equidistant surfaces reads as:
\begin{equation}\label{eq:eikonal}
\left\{ 
\begin{array}{l}
\left\|\nabla^{\epsilon}W(g)\right\|_{\epsilon}=1, \textrm{ if }g\neq e, \\
W(e)=0.
\end{array}
\right.
\end{equation} 
When $\epsilon\downarrow 0$ this system becomes the SR-eikonal system in \cite[eq.3]{bekkers} where it was proved that the unique viscosity solution is indeed the geodesic distance map from the origin $W(g)=d_\epsilon(g,e)$.
% viscosity solution $W(g)$ of the following boundary value problem \cite{bekkers}:
%\begin{equation}
%\left\{ 
%\begin{array}{l}
%\sqrt{\beta^{-2}|\mathcal{A}_{1}W(g)|^2 + |\mathcal{A}_{2}W(g)|^2}-1=0, \textrm{ if }g\neq e, \\
%W(e)=0.
%\end{array}
%\right.
%\end{equation} 
%The first equation of this system, a sub-Riemannian version of an eikonal equation, characterizes the propagation of geodesically equidistant surfaces $S_t$ (at distance $t$) as level sets of  $W(g)$:
%\begin{equation} \label{St}
%\mathcal{S}_{t}=\{g \in SE(2)\;|\; W(g)=t\}.
%\end{equation}
Then, SR-geodesics are the solutions $\gamma_b(t)$ of the following ODE system for backtracking:
\begin{equation}\label{eq:back}
\left\{\begin{array}{l}
\dot{\gamma_b}(t)= -\nabla^\epsilon (W(\gamma_b(t))), \,\, t\in[0,T]\\
 \gamma_b(0)=g.
\end{array}\right.
\end{equation}
%\begin{equation}\label{steepest}
%\begin{array}{ll}
%\dot{\gamma}(t)= - \frac{\left.\mathcal{A}_{1}W\right|_{\gamma(t)}}{(\beta \, \mathcal{C}(\gamma(t)))^2} \left.\mathcal{A}_{1}\right|_{\gamma(t)}-
%\frac{\left.\mathcal{A}_{2}W\right|_{\gamma(t)}}{
%(\mathcal{C}(\gamma(t)))^2} \left.\mathcal{A}_{2}\right|_{\gamma(t)}, & \gamma(0)=g.
%\end{array}
%\end{equation}
 %In \cite{bekkers} it was proved that the surfaces $S_t$ are indeed sub-Riemannian spheres of radius $t$, i.e.$S_t=\left\{g\in SE(2)| d(g,e)=t\right\}$. 
%\subsection{Riemannian approximation of the Sub-Riemannian problem}

\textbf{The metric matrix-representation in the Cartesian frame.}  A symmetric matrix $M_{\epsilon}$ representing the anisotropic metric in the frame $\{\partial x,\partial y,\partial\theta\}$ can be obtained by a basis transformation from the varying frame $\{ X_1,X_2,X_3\}$ (see \cite[Sec. 2.6]{cittisarti}):\begin{equation}
M_{\mathcal{\epsilon}} = \left(\begin{array}{rrr}
\cos\theta & 0 & -\sin\theta  \\
\sin\theta  & 0 & \cos\theta  \\
0& 1 & 0
\end{array}
\right) \left(\begin{array}{ccc}
\mathcal{C}^2\beta^2 & 0 & 0 \\
0& \mathcal{C}^2 &0 \\
0& 0 & \epsilon^{-2}\mathcal{C}^2\beta^2
\end{array}
\right) \left(\begin{array}{rrr}
\cos\theta & 0 & -\sin\theta  \\
\sin\theta  & 0 & \cos\theta  \\
0& 1 & 0
\end{array}
\right)^T .
\end{equation}
Here the diagonal matrix in the middle encodes the anisotropy between the $X_i$ directions while the other 2 rotation matrices are the basis transformation. Note that the columns are the coordinates of the varying frame in the fixed frame, e.g. $X_1=\cos\theta \partial x+\sin\theta \partial y+ 0\partial \theta $. Then, the metric tensor can be written as 
%\begin{equation}
$G_{\epsilon}(\dot\gamma,\dot\gamma)=\dot\gamma(t)M_{\epsilon} \dot\gamma(t) $
%\end{equation}
, with $\dot\gamma(t)=(\dot x(t),\dot y(t),\dot \theta (t))$.
Using this convention, the eikonal system (\ref{eq:eikonal}) in the fixed frame can be expressed as:
\begin{equation}\label{eq:eikfix}
\left\{\begin{array}{l}\nabla^T W(g)  M_{\epsilon}^{-1}\nabla W(g) = 1, \textrm{ if }g\neq e, \\
W(e)=0,
\end{array}\right.
%\, \mathrm{ with } \, \nabla W = \left( \frac{\partial W}{\partial x},\frac{\partial W}{\partial y}, \frac{\partial W}{\partial \theta} \right)^T
\end{equation}
 where $\nabla = (\partial_x,\partial_y,\partial_\theta)^T$ is the usual Euclidean gradient. The same holds for the backtracking equation (\ref{eq:back}) which writes:
\begin{equation}\label{eq:bacfix}
\left\{
\begin{array}{l}
\dot\gamma_b(t)= - M_{\epsilon}^{-1} \nabla W(\gamma_b(t)), \,\, t\in[0,T]\\
\gamma_b(0)=g.
\end{array}\right.
\end{equation}
Note that when approaching the SR-case, the $\lim\limits_{\epsilon\downarrow 0} M_{\epsilon}$ does not exist but the  $\lim\limits_{\epsilon\downarrow 0} M_{\epsilon}^{-1}$ is well defined in Eq.~(\ref{eq:bacfix}).\\ %It is also immediate to verify that $M_{1,1}$ is the 3D identity matrix when $\beta=1$ so in this case eqs. (\ref{eq:eikfix}) and (\ref{eq:bacfix}) are the usual Euclidean eikonal and backtracking systems respectively.\\
%\subsubsection{Anisotropic Fast Marching.}
%TEXT NEEDED HERE
%It is possible to construct a family of anisotropic Riemannian metric tensors as: $G_\epsilon=\beta^2 \omega^1 \otimes \omega^1 + \omega^2 \otimes \omega^2 + \beta^2 \epsilon^{-2} \omega^3 \otimes \omega^3  $, which bridges the SR-metric $G$ defined in eq. .. to the full Riemannian metric tensor.
%\begin{equation}
%G=\lim_{\epsilon\rightarrow 0} G_\epsilon
%\end{equation}

\textbf{Anisotropic Fast Marching.} We can now immediately identify Eq.~(\ref{eq:eikfix}) with \cite[eq.~0.1]{mirebeau} and then solve the eikonal system via the FM-LBR method. Our empirical tests show that $\epsilon=0.1$, which gives an anisotropy ratio $\kappa=0.01$ (see \cite[eq.~0.5]{mirebeau}), is already a good approximation of the SR-case and is the value used in the following experiments.

\section{Experiments and Applications} \vspace{-1.5ex}
%\subsubsection{Implementation.}
%TEXT NEEDED HERE

\textbf{Validation via comparison in uniform cost case.} The exact solutions of the SR-geodesic problem for the case $\mathcal{C}=1$ are known (see \cite{yuri} for optimal synthesis of the problem). Therefore, and similar to what is done in \cite{bekkers}, we consider this case as our golden standard.
Here, we want to compare both the computational time and the accuracy achieved in the calculation of the discrete SR-distance map $W(g)$, which is the solution of the eikonal system (\ref{eq:eikonal}) when $\epsilon \downarrow 0$.

\begin{figure}
\centering
\includegraphics[height=.28\textwidth]{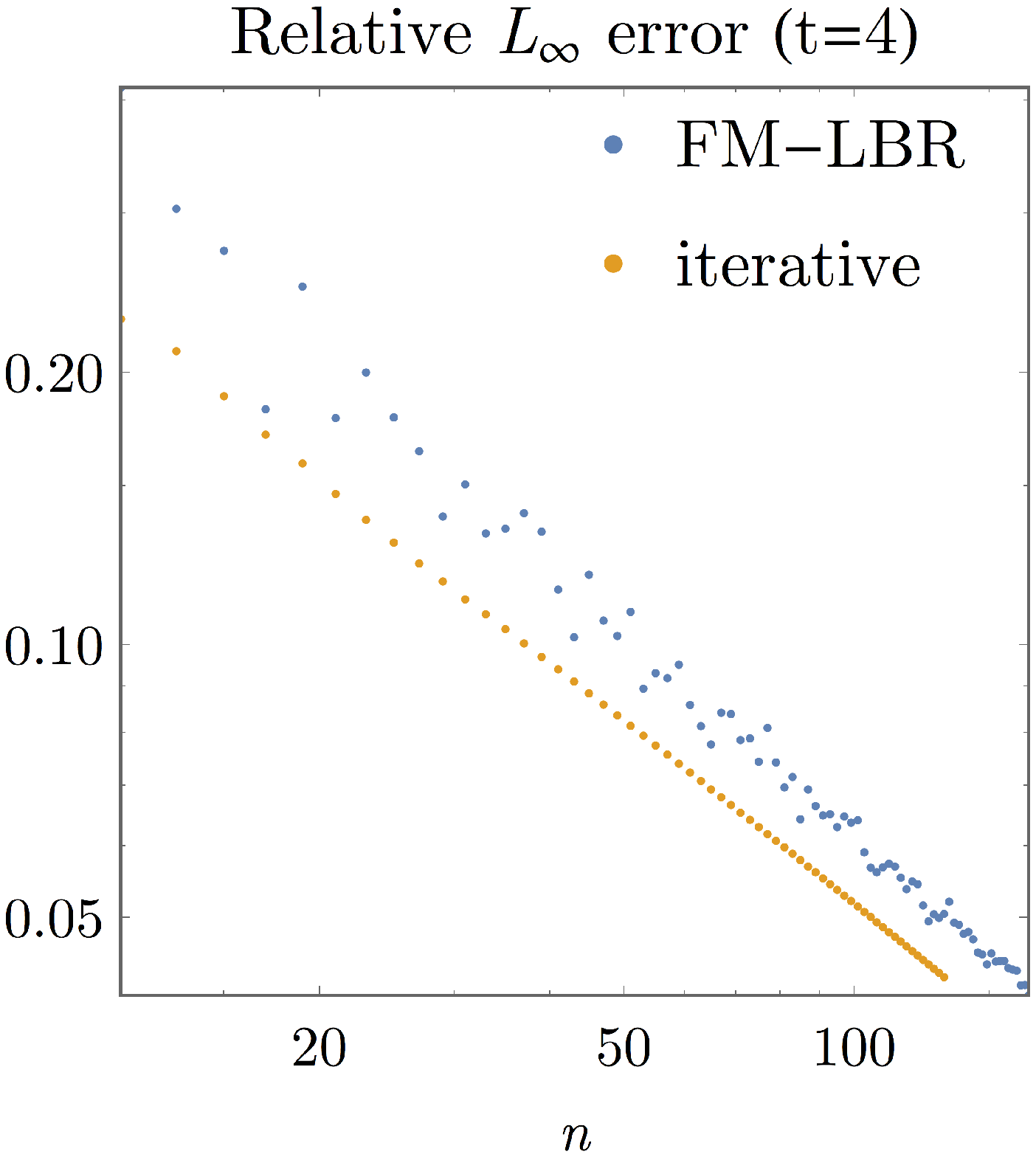}
\includegraphics[height=.28\textwidth]{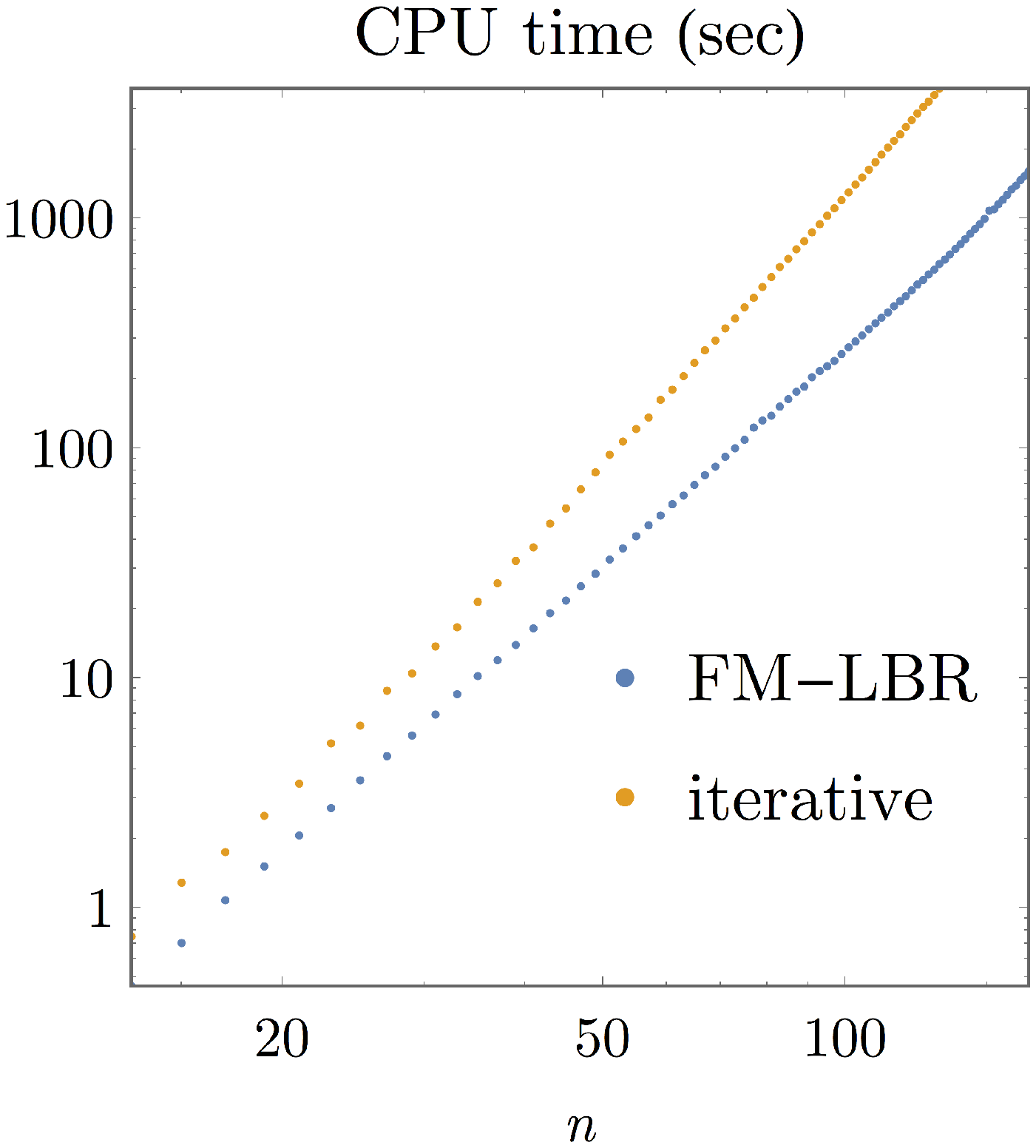}
\includegraphics[height=.3\textwidth]{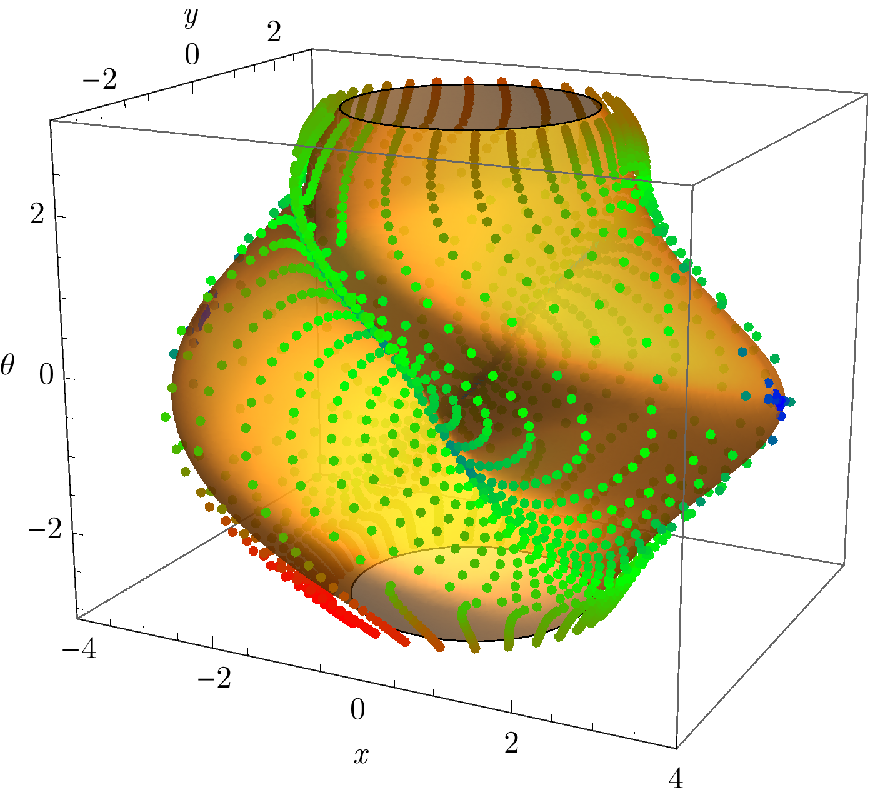}
\caption{Validation via comparison in the uniform cost case. The experiment (illustrated on the left, see the text) shows that even though the iterative method in \cite{bekkers} is more accurate we can still achieve with the FM-LBR method better results  and with less CPU effort, just by increasing the grid sampling.}
\label{fig:validation}
\end{figure}

Let us set $\beta=1$ and consider a grid $\mathcal{G}_s=\left\{(x_i,y_j,\theta_k) \in SE(2)\right\}$ with uniform step size $s$,  where $x_i=is$, $y_j=js,$ $\theta_k=k s,$ the indices $i,j,k\in\mathbb{N}$ such that $|x_i|\leq 2\pi,$ $|y_j|\leq 2\pi$ and $-\pi+s\leq\theta_k\leq \pi$. Then we compute the discrete geodesic distance map $W(g)$ on $\mathcal{G}_s$ using the iterative method in \cite{bekkers} and the FM-LBR. % as we have proposed in this paper. %, for a range of values of the spatial step size $s$. 
In order to measure the accuracy of the achieved solutions we follow the method explained in detail in \cite{bekkers}. There, by solving the initial value problem from the origin $e$, a set of arc length parametrized SR-geodesics is computed such that SR-spheres of radius $t$ are densely sampled. Then, the endpoints $g=(x,y,\theta)$ of each geodesic is stored in a list together with its length $t$.  Finally, we define the max relative error as $E_{\infty}(t)=\max(|W(g)-t|/t) $ where the max is taken over all the endpoints $g$ and where the value of $W(g)$ is obtained by bi-linearly interpolating the numerical solutions of eq. (\ref{eq:eikonal}) computed on the grid $\mathcal{G}_s$. The results and comparisons are presented in Fig.~\ref{fig:validation}. Here we solved the eikonal equation in increasingly finer grids  $\mathcal{G}_s$ obtained by setting the step size $s=\pi/n, n\in \mathbb{N}^+$. Note that the size of $\mathcal{G}_s$ is then $(2n+1)(2n+1)(n-1)$.  The graph in Fig.~\ref{fig:validation}(left) shows the comparison of the accuracy achieved in the computation of the SR-sphere of radius $t=4$ when $n$ increases. The behaviour for SR-spheres of different radius is similar. The CPU time is compared in Fig.~\ref{fig:validation}(center). The 3rd plot illustrates the method for computing the accuracy. The orange surface is the SR-sphere of radius $t=4$ computed with the FM-LBR method on a grid corresponding to $n=101$. The points are the geodesic endpoints, their color is proportional to the error of the FM-LBR (blue-low, green-medium, red-high error). The first observation is that even though the iterative method is more accurate, both methods seem to have the same order of convergence (the slope in the log-log graphs) when the grid resolution increases. This seems reasonable as both methods use first order approximations of the derivatives. Also, we hypothesise that the offset in favour of the iterative method is due to the Riemannian approximation of the SR-metric (i.e. selecting $\epsilon =0.1$), but this needs further investigation. The second key observation is that the CPU time increases dramatically with $n$ for the iterative method. Therefore, it is clear that we can achieve the same accuracy using the FM-LBR but with much less computational effort, which is of vital importance in the subsequent application. 

\begin{figure}
\centering
\includegraphics[width=.8\textwidth]{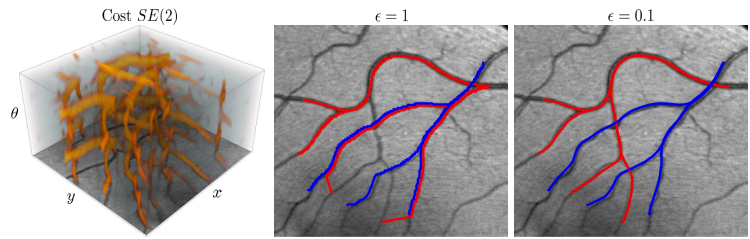}
\caption{Tracking of blood vessels in retinal images via cost adaptive SR-geodesics (see experiment details in \cite{bekkers}). \textbf{Left:} the cost obtained from the image (orange indicates locations with low cost). \textbf{Center:} Tracking in the full ($\epsilon=1$) Riemannian case.  \textbf{Right:}  Tracking with the approximated ($\epsilon=0.1$) SR-geodesics. }
\label{fig:retinal}
\end{figure}

\textbf{Application to retinal vessel tree extraction.} The analysis of the blood vessels in images of the retina  provides with early biomarkers of diseases such as diabetes, glaucoma or hypertensive retinopathy \cite{ikram}. For these studies, it is important to track the structural vessel tree, a difficult task especially because of the crossovers and bifurcations of the vessels. Some existing techniques \cite{Pechaud,bekkers} rely in considering an extended (orientation and/or scale) domain  to deal with this issue. Moreover, in \cite{bekkers} promising results were obtained  by formulating the vessel extraction as a $SE(2)$ SR-curve optimization problem with external cost obtained through some wavelet-like transformation of the 2D images. In the previous experiment, we have shown that our proposed FM-LBR based implementation computes in practice the same geodesics as the iterative method in \cite{bekkers}. Therefore, by simply replacing the eikonal equation solver in \cite{bekkers} we can obtain the same results but with a dramatic decrease of the computational demands (both CPU time and memory).  The example in Fig.~\ref{fig:retinal}  shows the advantages of considering the SR-metric in performing the vessel tree extraction. In this case (a patch of size 200x200, with 64 orientations considered) the iterative method computed the distance map in approximately 1 hour while the FM-LBR did the same in 20 seconds. For the experiments details we refer to \cite{bekkers}, for more examples see \url{www.bmia.bmt.tue.nl/people/RDuits/Bekkersexp.zip}. 
%\url{https://www.dropbox.com/s/yud42w9bgjcj5vs/patch_AFM.zip}.
%\begin{figure}
%\centering
%\includegraphics[width=.8\textwidth]{img/cost1.eps}\\
%\includegraphics[width=.8\textwidth]{img/example1.eps}
%\caption{Tracking of blood vessels in retinal images via cost adaptive SR-geodesics (see experiment details in \cite{bekkers}). \textbf{Top row:} the input image and the constructed cost (orange indicates locations with low cost). \textbf{Bottom row:} Better results are obtained when the SR-metric is taken into account. }
%\label{fig:retinal}
%\end{figure}
\section{Conclusions} \vspace{-1.5ex}
%Even though in this conference paper we only address the SR geometry in $SE(2)$ and the applications superficially, in our opinion it is a rich and powerful tool for image analysis related problems and a relatively newly developed field which is the mathematical (geometrical) modelling of the brain both at a neurophysiological level (modelling long range horizontal connectivity in V1, mean field equations that take into account this geometry) and its relation with the low level visual perception (association fields, figure ground segmentation taking into account self-occlusions).   Although some interesting application of SR geometry in the rototranslation group 
Over the last decade, some authors \cite{cittisarti,duits,franken} have shown the advantages of considering the roto-translation group embedded with a SR-geometry as a powerful, rich tool in some image analysis related problems or for the geometrical modelling of the visual perception. In our opinion, %mainly 3 
2 obstacles have prevented this framework to become more popular amongst engineers: %the high level of technicality inherent to the tool itself, 
the expensive computational demands involved (resulting of considering the extended orientation space) and the lack of efficient numerical methods able to deal with the extreme (degenerated) anisotropy imposed by the SR-metric. % especially when consider these non-linear of the Hamilton Jacobi equations depending on the external cost. 
%Regarding the 1st obstacle we have set up this paper using a concise notation that would allow readers with a background on PDE-based image processing to link the classical eikonal-based approach for computing geodesics with this SR extension to the $SE(2)$ moving frame.% embedded in the (sub-)Riemannian metric. 
%The 2nd and 3rd obstacle 
These obstacles are addressed by the main contribution of this work, which is solving (up to our knowledge for the first time) the SR-geodesic problem using a Fast Marching based implementation. To be able to achieve this, we rely on the FM-LBR solver recently introduce in \cite{mirebeau} and show that even when relaxing the SR-restriction by a Riemannian approximation of the metric we achieve excellent numerical convergence, but much faster than with the approach in \cite{bekkers}. Regarding the retinal imaging application our promising preliminary studies suggest that it is at least feasible to aim for a full vessel tree segmentation as the solution of a single optimization problem, but this requires further investigation. Future work will pursue extension of this method to the  3D-rototranslation group  $SE(3)$ and the applications in neuroimaging and 3D-vessel segmentation.  
{\begin{small}\end{small}}

%\newpage
%Extra (if there is space):

\end{document}